\documentclass{icmart}
\usepackage{amssymb, verbatim}

\usepackage{tikz}
\usetikzlibrary{arrows, shapes, calc}
\tikzstyle{every picture}+=[remember picture]

\newtheorem{theorem}{Theorem}

\newtheorem{conj}[theorem]{Conjecture}
\newtheorem{defn}[theorem]{Definition}

\theoremstyle{definition}

\newtheorem{example}[theorem]{Example}
\renewcommand{\[}{\begin{equation}}
\renewcommand{\]}{\end{equation}}
\newcommand{\ddb}{\sqrt{-1}\partial\overline{\partial}}
\newcommand{\Bl}{\mathrm{Bl}}

\def\XXint#1#2#3{{\setbox0=\hbox{$#1{#2#3}{\int}$ }
\vcenter{\hbox{$#2#3$ }}\kern-.6\wd0}}

\title{Extremal K\"ahler metrics}
\author{G\'abor Sz\'ekelyhidi}
\contact[gszekely@nd.edu]{Department of Mathematics, University of Notre Dame, Notre
  Dame, IN 46556}

\begin{document}

\begin{abstract}
  This paper is a survey of some recent progress on the study of Calabi's extremal
  K\"ahler metrics. We first discuss the Yau-Tian-Donaldson
  conjecture relating the existence of extremal metrics to an
  algebro-geometric stability notion and we give some example settings
  where this conjecture has been established. We then turn to the
  question of what one expects when no extremal metric exists.   
\end{abstract}

\maketitle

\section{Introduction}
A basic problem in differential geometry is to find ``best'' or
``canonical'' metrics on smooth manifolds. The most famous example is
the classical uniformization theorem, which says that every closed
2-dimensional manifold admits a metric with constant curvature, and
moreover this metric is essentially unique in its conformal
class. Calabi's introduction of extremal K\"ahler metrics~\cite{Cal82}
is an attempt at finding a higher dimensional
generalization of this result, in the setting of K\"ahler
geometry.

There are of course other ways in which one could attempt to
generalize the uniformization theorem to higher dimensional
manifolds. One possibility is the Yamabe problem in the
context of conformal geometry. This says~\cite{LP87} that on a closed
manifold of arbitrary dimension, every conformal class admits a metric
of constant scalar curvature. Moreover this metric is often, but not
always, unique up to scaling. A different generalization to
3-dimensional manifolds is given by Thurston's geometrization
conjecture, established by Perelman~\cite{Per02}. In this case the
goal is to find metrics of constant curvature on a 3-manifold, but this is too
ambitious. Instead it turns out that every 3-manifold can be
decomposed into pieces each of which admits one of 8 model
geometries.

The search for extremal K\"ahler metrics can be thought of as a
complex analogue of the Yamabe problem, where we try to find canonical
representatives of a given K\"ahler class, rather than a conformal
class. In both cases the effect of restricting the space of metrics
that we allow results in the problems reducing to scalar equations
involving the conformal factor, and the K\"ahler potential,
respectively. We will see, however, that in contrast with the Yamabe
problem extremal metrics do not always exist, and in these cases one
can hope to find a canonical ``decomposition'' of the manifold into
pieces somewhat reminiscent of the geometrization of 3-manifolds.

In order to define extremal metrics, let $M$ be a compact complex
manifold of dimension $n$, equipped with a K\"ahler class $\Omega\in
H^2(M,\mathbf{R})$. Denote by $\mathcal{K}_\Omega$ the set of K\"ahler
metrics in the class $\Omega$.
\begin{defn}
  An extremal metric is a critical point of the
Calabi functional
\[ \begin{aligned}
  \mathrm{Cal} : \mathcal{K}_\Omega &\to \mathbf{R} \\
  \omega &\mapsto \int_M (S(\omega) -
  \underline{S})^2\,\omega^n,
  \end{aligned}\]
where $S(\omega)$ is the scalar curvature of $\omega$, and
$\underline{S}$ is the average of $S(\omega)$ with respect to the
volume form $\omega^n$. Note that $\underline{S}$ is independent of
the choice of $\omega\in \mathcal{K}_\Omega$.
\end{defn}

Calabi~\cite{Cal82} has shown that $\omega$ is an extremal metric if
and only if the gradient $\nabla S(\omega)$ is a holomorphic vector
field. Since most complex manifolds
do not admit any non-trivial holomorphic vector fields, most extremal
metrics are constant scalar curvature K\"ahler (cscK) metrics. A
particularly important special case is when the first Chern class
$c_1(M)$ is proportional to the K\"ahler class $\Omega$. If $c_1(M) =
\lambda\Omega$, and $\omega\in\mathcal{K}_\Omega$ is a cscK metric,
then it follows that
\[ \mathrm{Ric}(\omega) = \lambda\omega, \]
and so $\omega$ is a K\"ahler-Einstein metric. 

It is known that any two extremal metrics in a fixed K\"ahler class
are isometric (see Chen-Tian~\cite{CT05_1}), which makes extremal
metrics good candidates for being canonical metrics on K\"ahler
manifolds. On the other hand, not every K\"ahler class admits an
extremal metric, the first examples going back to
Levine~\cite{Levine85} of manifolds which do not admit extremal
metrics in any K\"ahler class. The basic problems are therefore to understand
which K\"ahler classes admit extremal metrics, and what we can say
when no extremal metric exists.

The most interesting case of the existence question is when $\Omega=c_1(L)$
is the first Chern class of a line bundle, and consequently $M$ is a projective
manifold. In this case the Yau-Tian-Donaldson conjecture predicts that
the existence of an extremal metric is related to the stability of the
pair $(M,L)$ in the sense of geometric invariant theory. In
Section~\ref{sec:Kstab} we will discuss two such notions of
stability: K-stability, and a slight refinement of it which we call
$\widehat{K}$-stability. 

As a consequence of work of Tian~\cite{Tian97}, Donaldson~\cite{Don01,
  Don05}, Mabuchi~\cite{Mab08}, Stoppa~\cite{Sto08}, Stoppa-Sz\'ekelyhidi~\cite{SSz09},
Paul~\cite{Paul12}, Berman~\cite{Ber12}, and others, there are now
many satisfactory results that show that the existence of an extremal metric
implies various notions of stability. The converse direction, however,
is largely open. In Section~\ref{sec:existence} we will discuss two
results in this direction. One is the recent breakthrough of
Chen-Donaldson-Sun~\cite{CDS12} on K\"ahler-Einstein metrics with
positive curvature, and the other is work of the author on extremal
metrics on blowups. 

Finally in Section~\ref{sec:noextremal} we turn to what is to be expected when no extremal metric
exists, i.e. when a pair $(M,L)$ is unstable.
It is still a natural problem to try minimizing the Calabi
functional in a K\"ahler class, and we will discuss a conjecture due to
Donaldson relating this to finding the optimal way to destabilize
$(M,L)$. We will give an example where this can be interpreted as the
canonical decomposition of the manifold alluded to above. 

\subsection*{Acknowledgements}
Over the years I have benefited from conversations about extremal
metrics with many people. In particular I would like to thank Simon
Donaldson, Duong Phong, Julius Ross, Jacopo
Stoppa, Richard Thomas and Valentino Tosatti for many useful
discussions. The work presented in this survey was partially supported
by the NSF. 

\section{The Yau-Tian-Donaldson conjecture}\label{sec:Kstab}
 It is a conjecture going back
to Yau (see e.g. \cite{Yau93}) that if $M$ is a Fano manifold,
i.e. the anticanonical line bundle $K_M^{-1}$ is ample, then $M$
admits a K\"ahler-Einstein metric if and only if $M$ is stable in the
sense of geometric invariant theory. Tian~\cite{Tian97} introduced the
notion of K-stability as a precise candidate of such a stability
condition and showed that it is necessary for the existence of a
K\"ahler-Einstein metric. Donaldson~\cite{Don97, Don02} generalized
the conjecture to pairs $(M,L)$ where $L\to M$ is an ample line
bundle, not necessarily equal to the anticanonical bundle. More
precisely, Donaldson formulated a more algebraic version of
K-stability, and conjectured that K-stability of the pair $(M,L)$ is
equivalent to the existence of a cscK metric in the class $c_1(L)$. We
start by giving a definition of Donaldson's version of K-stability.

\begin{defn} A \emph{test-configuration} for $(M,L)$  of exponent $r$
  is a $\mathbf{C}^*$-equivariant, flat, polarized family
  $(\mathcal{M},\mathcal{L})$ over $\mathbf{C}$, with generic fiber
  isomorphic to $(M, L^r)$. 
\end{defn}

The central fiber $(M_0, L_0)$ of a test-configuration 
has an induced $\mathbf{C}^*$-action, and we write $A_k$ for the
infinitesimal generator of this action on $H^0(M_0, L_0^k)$. In other
words, the eigenvalues of $A_k$ are the weights of the action. There
are expansions
\[ \begin{aligned}
  \dim H^0(M_0, L_0^k) &= a_0k^n + a_1k^{n-1} + O(k^{n-2}) \\
  \mathrm{Tr}(A_k) &= b_0k^{n+1} + b_1k^n + O(k^{n-1}) \\
  \mathrm{Tr}(A_k^2) &= c_0k^{n+2} + O(k^{n+1}).
\end{aligned}\]

\begin{defn}
 Given a test-configuration $\chi$ of exponent $r$ for $(M,L)$ as above, its Futaki
 invariant is defined to be
\[ \mathrm{Fut}(\chi) = \frac{a_1b_0 - a_0b_1}{a_0^2}. \]
The norm of $\chi$ is defined by
\[ \Vert\chi\Vert^2 = r^{-n-2}\left(c_0 - \frac{b_0^2}{a_0}\right), \]
where the factor involving $r$ is used to make the norm unchanged if
we replace $L$ by a power. 
\end{defn}

With these preliminaries, we can give a definition of K-stability.
\begin{defn}\label{defn:Kstab}
The pair $(M,L)$ is K-stable if $\mathrm{Fut}(\chi) > 0$ for all
test-configurations $\chi$ with $\Vert \chi\Vert > 0$. 
\end{defn}

The condition $\Vert\chi\Vert > 0$ is required to rule out certain
``trivial'' test-configurations. An alternative definition by
Li-Xu~\cite{LX11} requires $\mathcal{M}$ to be a normal variety
distinct from the product $M\times\mathbf{C}$, but the condition using
the norm $\Vert\chi\Vert$ will be more natural below, when discussing
filtrations.

The central conjecture in the field is the following.
\begin{conj}[Yau-Tian-Donaldson]\label{conj:YTD}
Suppose that $M$ has
no non-zero holomorphic vector fields. Then $M$ admits a cscK metric in
$c_1(L)$ if and only of $(M,L)$ is K-stable. 
\end{conj}

When $M$ admits holomorphic vector fields which can be lifted to $L$, then it is
never K-stable according to the previous definition, since in that case one can
find test-configurations with total space $\mathcal{M} =
M\times\mathbf{C}$, with a non-trivial $\mathbf{C}^*$-action whose
Futaki invariant is non-positive. In this case a variant of K-stability, called
K-polystability, is used, which rules out such ``product
test-configurations'' and is conjecturally equivalent to the existence
of a cscK metric, even when $M$ admits holomorphic vector fields. 
A further variant of K-stability, called relative K-stability, was defined by
the author~\cite{GSz04}, and  it is conjecturally related to the
existence of extremal metrics. In relative K-stability one only
considers test-configurations which are orthogonal to a maximal torus
of automorphisms of $M$ in a suitable sense. 

\begin{example}
  Let $(M,L) = (\mathbf{P}^1, \mathcal{O}(1))$. The family of conics
  $xz=ty^2$ for $t\in \mathbf{C}$ gives a test-configuration $\chi$ for
  $(M,L)$ of exponent 2, degenerating a smooth conic into the union of
  two lines (see Figure~\ref{fig:tc}). 
A small computation gives $\mathrm{Fut}(\chi)
  = 1/8$. It is not surprising that this is positive, since the
  Fubini-Study metric on  $\mathbf{P}^1$ has constant scalar
  curvature. 
  \begin{figure}[h]
    \centering
    \begin{tikzpicture}
      \draw[line width=1.5pt] (0,2) to [out=290, in=80] (-0.2,-0.7);
      \draw (0,-0.8) node[below] {$xz=ty^2$};
      \draw[style=dashed, ->]  (1,0.5) -- (3,0.5);
      \draw (2,0.5) node[below] {$t\to 0$};
      \draw[line width=1.5pt] (4,1.8) to [out=320, in=110] (6,0);
      \draw[line width=1.5pt] (4.6, -0.8) to [out=80, in=225] (6,1);
      \draw (5.5,-0.6) node[below] {$xz=0$};
    \end{tikzpicture}
\caption{A test-configuration degenerating a conic into two lines.}
\label{fig:tc}
    \end{figure}
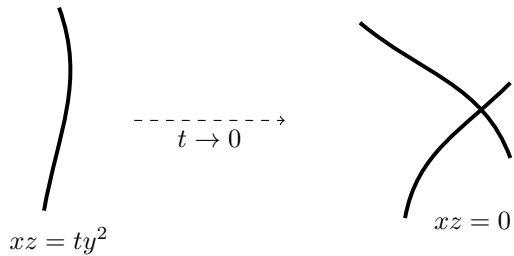
\end{example}

Calculations in
Apostolov-Calderbank-Gauduchon-T\o{}nnesen-Friedman~\cite{ACGT3}
suggest that K-stability might not be sufficient to ensure the existence of
a cscK metric in general. Indeed they construct examples where the
existence of an extremal metric is equivalent to the positivity of a
certain function $F$ on an interval $(a,b)$,
 while relative K-stability only ensures that $F$ is positive at rational
points $(a,b)\cap \mathbf{Q}$. It is thus natural to try to work with
a completion of the space of test-configurations in a suitable sense
in order to detect when this function $F$ vanishes at an irrational
point. This motivates the the author's work~\cite{GSz13} on
filtrations. 

\begin{defn}
  Let $R = \bigoplus_{k\geqslant 0} H^0(M,L^k)$ denote the homogeneous
  coordinate ring of $(M,L)$. A filtration of $R$ is a family of
  subspaces
\[ \mathbf{C} = F_0R \subset F_1R \subset\ldots \subset R, \]
 satisfying
\begin{enumerate}
\item $(F_iR)(F_jR) \subset F_{i+j}R$,
\item If $s\in F_iR$, and $s$ has degree $k$ piece $s_k\in H^0(M,L^k)$, then $s_k\in
  F_iR$. 
\item $R = \bigcup_{i\geqslant 0} F_i R$. 
\end{enumerate}
\end{defn}

Witt Nystr\"om~\cite{Ny10} showed that every test-configuration for
$(M,L)$ gives rise to a filtration. In fact the Rees
algebra of the filtration is the coordinate ring of the total space of
the test-configuration. On the other hand given any filtration $\chi$
of the homogeneous coordinate ring of $(M,L)$, we obtain a flag of
subspaces
\[ \{0\} = F_0R_r \subset F_1R_r \subset\ldots \subset R_r \]
of the degree $r$ piece $R_r = H^0(M,L^r)$ for all $r > 0$. In turn,
such a flag gives rise to a test-configuration of exponent $r$ for
$(M,L)$ -- by embedding $M$ into projective space using a basis of $H^0(M,L^r)$, and acting by  a
$\mathbf{C}^*$-action whose weight filtration is given by our flag. 
Therefore any filtration $\chi$ induces a sequence of
test-configurations $\chi^{(r)}$, where $\chi^{(r)}$ has exponent
$r$. It is natural to think of $\chi$ as the limit of the
$\chi^{(r)}$, and thus to define
\[\begin{aligned}
 \mathrm{Fut}(\chi) &= \liminf_{r\to\infty} \mathrm{Fut}(\chi^{(r)}) \\
 \Vert\chi\Vert &= \lim_{r\to\infty} \Vert\chi^{(r)}\Vert, 
\end{aligned} \]
where the limit can be shown to exist. The main difference between
filtrations arising from test-configurations, and general filtrations,
is that the Rees algebras of the latter need not be finitely
generated. 

In terms of filtrations we define the following stability notion,
which is stronger than K-stability.
\begin{defn}
The pair $(M,L)$ is $\widehat{K}$-stable, if $\mathrm{Fut}(\chi) > 0$
for all filtrations of the homogeneous coordinate ring of $(M,L)$
satisfying $\Vert\chi\Vert > 0$. 
\end{defn}

In view of the examples of Apostolov-et. al. that we mentioned above,
it may be that in the Yau-Tian-Donaldson conjecture one should
assume $\widehat{K}$-stability instead of K-stability. One direction
of this modified conjecture 
has been established by Boucksom and the author~\cite{GSz13}. 
\begin{theorem}\label{thm:Khatstable}
Suppose that $M$ has no non-zero holomorphic vector fields. If $M$
admits a cscK metric in $c_1(L)$, then $(M, L)$ is
$\widehat{K}$-stable. 
\end{theorem}
The analogous result for K-stability was shown in
Stoppa~\cite{Sto08}, building on work of Donaldson~\cite{Don05} which
we will see in Theorem~\ref{thm:Callower}, and
Arezzo-Pacard~\cite{AP06} which we will discuss in
Section~\ref{sec:existence}. In the proof of
Theorem~\ref{thm:Khatstable} the main additional ingredient is the use of
the Okounkov body~\cite{Ok96, BC09, LM09, Ny10}. 
Note that when $L=K_M^{-1}$, then
related results were shown by Tian~\cite{Tian97} and
Paul-Tian~\cite{PT06}. It is likely that a result analogous to
Theorem~\ref{thm:Khatstable} can be shown for extremal metrics along
the lines of \cite{SSz09}.

\section{Some existence results}\label{sec:existence}
In this section we discuss two special cases, where the
Yau-Tian-Donaldson conjecture has been verified. 

\subsection*{K\"ahler-Einstein metrics}
We first focus on K\"ahler-Einstein metrics, i.e. when $c_1(M)$ is
proportional to $c_1(L)$. When $c_1(M) = 0$, or $c_1(M) < 0$, then the
celebrated work of Yau~\cite{Yau78} implies that $M$ admits a
K\"ahler-Einstein metric, and a stability condition does not need to
be assumed (see also Aubin~\cite{Aub78} for the case when $c_1(M) <
0$).

In the remaining case, when $c_1(M) > 0$, i.e. $M$ is Fano, it was known from
early on (see e.g. Matsushima~\cite{Mat57}) that a K\"ahler-Einstein
metric does not always exist, and Yau conjectured that the existence
is related to stability of $M$ in the sense of geometric invariant
theory.
 Tian~\cite{Tian90} found all two-dimensional $M$ which admit a
 K\"ahler-Einstein metric, and in \cite{Tian97} he formulated
the notion of K-stability, which he conjectured to be equivalent to
the existence of a K\"ahler-Einstein metric.  The main difference
between Tian's notion of K-stability and the one in
Definition~\ref{defn:Kstab} is that Tian's version of K-stability only
requires $\mathrm{Fut}(\chi)> 0$ for very special types of
test-configurations with only mild singularities. In particular their
Futaki invariants can be computed differential geometrically using the
formula Futaki~\cite{Fut83} originally used to define his
invariant. By the work of Li-Xu~\cite{LX11} it turns out that in the
Fano case
Tian's notion of K-stability is equivalent to the a priori
stronger condition of Definition~\ref{defn:Kstab}. 

Recently, Chen-Donaldson-Sun~\cite{CDS12, CDS13_1, CDS13_2, CDS13_3}
have proved Conjecture~\ref{conj:YTD} for Fano manifolds:
\begin{theorem}
  Suppose that $M$ is a Fano manifolds and $(M, K_M^{-1})$ is
  K-polystable. Then $M$ admits a K\"ahler-Einstein metric.
\end{theorem}

To construct a K\"ahler-Einstein metric, the continuity method is
used, with a family of equations of the form
\[ \label{eq:conecont} \mathrm{Ric}(\omega_t) = t\omega_t + \frac{1-t}{m}[D], \]
where $D$ is a smooth divisor in the linear system $|mK_M^{-1}|$, and
$[D]$ denotes the current of integration. More precisely, a metric
$\omega_t$ is a solution of $\eqref{eq:conecont}$ if
$\mathrm{Ric}(\omega_t)=t\omega_t$ on $M\setminus D$, while $\omega_t$
has conical singularities along $D$ with cone angle
$\frac{2\pi}{m}(1-t)$. One then shows that Equation~\ref{eq:conecont}
can be solved for $t\in [t_0, T)$ for some $t_0, T > 0$ (see
Donaldson~\cite{Don11_1} for the openness statement), and the 
question is what happens when $t\to T$.

One of the main results of the
work of Chen-Donaldson-Sun is, roughly speaking, that along a
subsequence the manifolds $(M, \omega_{t_i})$ have a 
Gromov-Hausdorff limit $W$ which is a $\mathbf{Q}$-Fano variety, such
that the divisor $D\subset M$ converges to a Weil divisor $\Delta$,
and $W$ admits a weak K\"ahler-Einstein metric with conical
singularities along $\Delta$ (defined in an appropriate
sense). Moreover, there are embeddings $\phi_i:
M\to\mathbf{P}^N$ and
$\phi: W\to\mathbf{P}^N$
into a sufficiently large projective space, such that the pairs $(\phi_i(M), \phi_i(D))$
converge to $(\phi(W), \phi(\Delta))$ in an algebro-geometric
sense. In this case either $(\phi(W), \phi(\Delta))$ is in the
$SL(N+1,\mathbf{C})$-orbit of the $(\phi_i(M), \phi_i(D))$ in which
case we can solve Equation~\eqref{eq:conecont} for $t=T$, or
otherwise we can find a test-configuration for $(M,D)$ with central
fiber
$(W, \Delta)$ to show that $(M,K_M^{-1})$ is not K-stable. Note that here one
needs to extend the theory described in Section~\ref{sec:Kstab} to
pairs $(M,D)$ resulting in the notion of log
K-stability~\cite{Don11_1}. 

The fact that a sequence of solutions to Equation~\eqref{eq:conecont}
has a Gromov-Hausdorff limit which is a $\mathbf{Q}$-Fano variety
originates in work of Tian~\cite{Tian90} on the 2-dimensional case,
and it is essentially equivalent to what Tian calls the ``partial $C^0$-estimate''
being satisfied by such a sequence of solutions. This partial
$C^0$-estimate was first shown in dimensions greater than 2 by
Donaldson-Sun~\cite{DS12} for sequences of K\"ahler-Einstein metrics,
and their method has since been generalized to many other settings:
Chen-Donaldson-Sun~\cite{CDS13_2, CDS13_3} to solutions of
\eqref{eq:conecont}; Phong-Song-Sturm~\cite{PSS12} for sequences of K\"ahler-Ricci
solitons; Tian-Zhang~\cite{TZ13} along the K\"ahler-Ricci flow in
dimensions at most 3; the author~\cite{Sz13_1} along Aubin's continuity
method; Jiang~\cite{Jiang13} using only a lower bound for the Ricci
curvature, in dimensions at most 3. Note that Tian's original
conjecture on the partial $C^0$-estimate is still open in dimensions
greater than 3 -- namely we do not yet understand Gromov-Hausdorff limits of Fano
manifolds under the assumption of only a positive lower bound on the Ricci
curvature.

To close this subsection we mention a possible further result along
the lines of Chen-Donaldson-Sun's work. As we described above, if $M$ does
not admit a K\"ahler-Einstein metric, then a sequence of solutions to
Equation~\eqref{eq:conecont} will converge to a weak conical
K\"ahler-Einstein metric on a pair $(W,\Delta)$ as $t\to T$. Suppose
$T < 1$. We can think of this metric as a suitable weak solution to
the equation
\[ \mathrm{Ric}(\omega_t) = t\omega_t + \frac{(1-t)}{m}[\Delta] \]
for $t=T$ on the space $W$. Since the pair $(W,\Delta)$ necessarily has a
non-trivial automorphism group, we cannot expect to solve this
equation for $t > T$, however it is reasonable to expect that we can
still find weak conical K\"ahler-Ricci solitons, i.e. we can solve
\[ \label{eq:coneKRS}
\mathrm{Ric}(\omega_t) + L_{X_t}\omega_t = t\omega_t +
\frac{(1-t)}{m}[\Delta], \]
for some range of values $t > T$, with
suitable vector fields $X_t$ fixing $\Delta$. An extension of
Chen-Donaldson-Sun's work to K\"ahler-Ricci solitons, generalizing
Phong-Song-Sturm~\cite{PSS12}, could then be used to extract a limit
$(W_1, \Delta_1)$ as $t\to T_1$, with yet another (weak) conical
K\"ahler-Ricci soliton, and so on. Based on this heuristic argument we
make the following conjecture.

\begin{conj} We can solve Equation~\eqref{eq:coneKRS} up to $t=1$ by
  passing through finitely many singular times, changing the pair
  $(W,\Delta)$ each time. At $t=1$ we obtain a $\mathbf{Q}$-Fano
  variety $W_k$ admitting a weak K\"ahler-Ricci soliton. Moreover,
  there is a test-configuration for $(M,K_M^{-1})$ with central fiber
  $W_k$. 
\end{conj}

A further natural expectation would be that the K\"ahler-Ricci soliton
obtained in this way is related to the limiting behavior of the
K\"ahler-Ricci flow on $M$. Indeed, according to the
Hamilton-Tian conjecture (see
Tian~\cite{Tian97}), the K\"ahler-Ricci flow is expected to converge to a
K\"ahler-Ricci soliton with mild singularities. 

\subsection*{Blow-ups}
Beyond the K\"ahler-Einstein case there are very few general existence
results for cscK or extremal metrics. One example is the case of toric
surfaces, where Conjecture~\ref{conj:YTD} has been established by
Donaldson~\cite{Don08_1}, with an extension to extremal metrics by
Chen-Li-Sheng~\cite{CLS10}. In this section we will discuss a 
perturbative existence result for cscK metrics on blow-ups. 

Suppose that $\omega$ is a cscK
metric on a compact K\"ahler manifold $M$, and choose a point $p\in
M$. For all sufficiently small $\epsilon > 0$ the class
\[ \Omega_\epsilon = \pi^*[\omega] - \epsilon^2 [E] \]
is a K\"ahler class on the blowup $\Bl_p M$, where $\pi: \Bl_p M\to M$
is the blowdown map, and $[E]$ denotes the Poincar\'e dual of the
exceptional divisor. A basic question, going back to work of
LeBrun-Singer~\cite{LS94}, is whether $\Bl_p M$ admits a cscK (or
extremal) metric in the class $\Omega_\epsilon$ for sufficiently small
$\epsilon$. This problem was studied extensively by
Arezzo-Pacard~\cite{AP06, AP09}, and
Arezzo-Pacard-Singer~\cite{APS06}. See also Pacard~\cite{Pac10} for a
survey. The following is the most basic result in this direction.

\begin{theorem}[Arezzo-Pacard~\cite{AP06}]\label{thm:AP}
  Suppose that $M$ admits a cscK metric $\omega$, and it 
admits no non-zero holomorphic vector fields. Then there is an
$\epsilon_0 > 0$ such that $\Bl_pM$ admits a cscK metric in the class
$\Omega_\epsilon$ for all $\epsilon\in(0,\epsilon_0)$. 
\end{theorem}

\begin{figure}
\centering
      \begin{tikzpicture}
        \draw[line width=1.5pt] (0,0) to [out=90, in=180] (2,1)
        to [out=0, in=140] (6,0.5)
        to [out=320, in=110] (6.2,0.2);
        \draw[line width=1.5pt] (6.2,0.2) to [out=290, in=180] (6.35, 0.08)
        to [out=0, in=180] (6.6, 0.17)
        to [out=0, in=90] (6.8, 0)
        to [out=270, in=0] (6.6, -0.17)
        to [out=180, in=0] (6.35, -0.08)
        to [out=180, in=70] (6.2, -0.2);
        \draw[line width=1.5pt] (6.2, -0.2) to [out=250, in=0] (3,-2)
        to [out=180, in=270] (0,0);
        \draw[line width=1.5pt] (2,-0.7) to [out=330, in=225] (4,-0.5);
        \draw[line width=1.5pt] (2.3, -0.8) to [out=45, in=135] (3.6, -0.7);

        \draw (6.2, 0.2) to (8.2, 2);
        \draw (6.2, -0.2) to (8.2, -2);

       \draw[line width=1.5pt] (8,2.7) to [out=290, in=180] (9.7, 0.2)
        to [out=0, in=180] (10.5, 0.5)
        to [out=0, in=90] (10.8, 0)
        to [out=270, in=0] (10.5, -0.5)
        to [out=180, in=0] (9.7, -0.2)
        to [out=180, in=70] (8, -2.7);

        \draw[densely dotted] (8.2,2) to [out=260, in=100] (8.2,-2);
        \draw[densely dotted] (6.2,0.2) to [out=260, in=100]
        (6.2,-0.2);

        \draw (9.6, -2.7) node[below, text width=4cm] {scalar flat  metric $\eta$
                   (Burns-Simanca)};
        \draw (3,2) node {$\mathrm{Bl}_pM$};
        \draw (9.7, 1.6) node {$\mathrm{Bl}_0\mathbf{C}^n$};
       \draw (3, -2.2) node[below] {$\omega_\epsilon = \omega$ glued to
          $\epsilon^2\eta$};
      \end{tikzpicture}
   \caption{The construction of the approximate metric
     $\omega_\epsilon$.}
\label{fig:gluing}
 \end{figure}
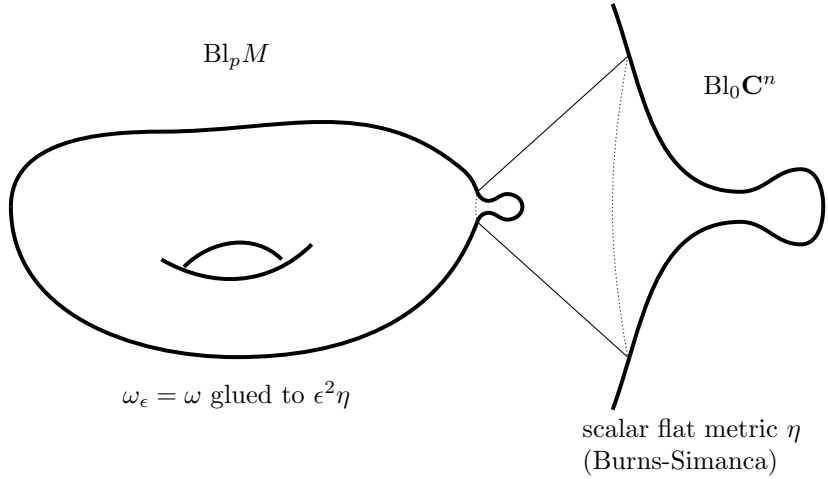

This result not only provides many new examples of cscK metrics, but
it is also a key ingredient in the proofs of results such as
Theorem~\ref{thm:Khatstable}. The construction of cscK metrics on
blowups is a typical example of a gluing theorem in geometric
analysis. First, one obtains a metric $\omega_\epsilon\in
\Omega_\epsilon$ on the blowup, by gluing the metric $\omega$ to a
scaled down version $\epsilon^2\eta$ of the scalar flat
Burns-Simanca~\cite{Sim91}
metric $\eta$ on $\Bl_0\mathbf{C}^n$. This is shown in
Figure~\ref{fig:gluing}. In a suitable weighted H\"older space the metric
$\omega_\epsilon$ is sufficiently close to having constant scalar
curvature, that one can perturb it to a cscK metric using a
contraction mapping argument, for sufficiently small $\epsilon$. 

When $M$ admits non-zero holomorphic vector fields, then the problem
becomes more subtle, since then $\Bl_pM$ may not admit a cscK (or even
extremal) metric for every point $p$. The problem was addressed by
Arezzo-Pacard~\cite{AP09} and Arezzo-Pacard-Singer~\cite{APS06} in the
case of extremal metrics, as well as the author~\cite{GSz10,
  GSz13_1}. For the case of cscK metrics the sharpest result from
\cite{GSz13_1} is as follows, showing that the Yau-Tian-Donaldson
conjecture holds for the pair $(\Bl_pM, \Omega_\epsilon)$ for sufficiently
small $\epsilon$. 

\begin{theorem}\label{thm:stabblowup} Suppose that $\dim M> 2$, $M$ 
  admits a cscK metric $\omega$, and
  $p\in M$. Then for sufficiently small $\epsilon > 0$, the blowup
  $\Bl_pM$ admits a cscK metric in the class $\Omega_\epsilon$ 
  if and only if $(\Bl_pM,\Omega_\epsilon)$ is K-polystable. 
\end{theorem}

For K-polystability to be defined algebraically, 
the class $\Omega_\epsilon$ should
be rational, but in fact a very weak version of K-polystability, which
can be defined for K\"ahler manifolds, is sufficient in this
theorem. Indeed what we can prove is that if $\Bl_pM$ does not admit a
cscK metric in the class $\Omega_\epsilon$ for sufficiently small
$\epsilon$, then there is a $\mathbf{C}^*$-action $\lambda$ on $M$
such that if we let 
\[q = \lim_{t\to 0} \lambda(t)\cdot p, \]
then the $\mathbf{C}^*$-action on $\Bl_q M$ induced by $\lambda$ has
non-positive Futaki invariant. In other words when $\epsilon$ is
sufficiently small, then it is enough to consider test-configurations
for $\Bl_pM$ which arise from one-parameter subgroups in the
automorphism group of $M$. While there are also existence results for
cscK metrics when $\dim M=2$, and also for general extremal metrics,
in these cases the precise relation with (relative) K-stability has
not been established yet. 

In the remainder of this section we will give a rough idea of the proof of
these existence results. The basic ingredient is the existence of a
scalar flat, asymptotically flat metric $\eta$ on $\Bl_0\mathbf{C}^n$
due to Burns-Simanca~\cite{Sim91}, of the form
\[ \eta = \ddb\Big[ |w|^2 + \psi(w)\Big] \]
on $\mathbf{C}^n\setminus\{0\}$, where
\[ \psi(w) = -|w|^{4-2n} + O(|w|^{2-2n}),\quad\text{as $|w|\to\infty$} \]
for $n>2$. Under the change of variables $w=\epsilon^{-1}z$, we have
\[\label{eq:e2eta}  \epsilon^2\eta = \ddb\Big[ |z|^2 +
\epsilon^2\psi(\epsilon^{-1}z)\Big]. \]
At the same time there are local coordinates near $p\in M$
for which the metric $\omega$ is of the form
\[\label{eq:omegaexpand} \omega = \ddb\Big[ |z|^2 + \phi(z)\Big], \]
where $\phi(z) = O(|z|^4)$. One can then use cutoff functions to glue
the metrics $\omega$ and $\epsilon^2\eta$ on the level of K\"ahler
potentials on the annular region $r_\epsilon < |z| < 2r_\epsilon$ for
some small radius $r_\epsilon$. The result is a metric $\omega_\epsilon \in
\Omega_\epsilon$ on $\Bl_pM$, which in a suitable weighted H\"older
space is very close to having constant scalar curvature if $\epsilon$
is small. It is important here that $\eta$ is scalar flat, since if it
were not, then $\epsilon^2\eta$ would have very large scalar curvature
once $\epsilon$ is small. 

When $M$ has no holomorphic vector fields, then one can show that for
sufficiently small $\epsilon$ this metric $\omega_\epsilon$ can be
perturbed to a cscK metric in its K\"ahler class, and this proves
Theorem~\ref{thm:AP}. Analytically the main ingredient in this proof
is to show that the linearization of the scalar curvature operator is
invertible, and to control the norm of its inverse in
suitable Banach spaces as $\epsilon\to 0$. 

The difficulty when $M$ has holomorphic vector
fields, or more precisely when when the Hamiltonian isometry group $G$ of
$(M,\omega)$ is non-trivial, is that the linearized operator will no
longer be surjective, since its cokernel can be identified with the
Lie algebra $\mathfrak{g}$ of $G$. One way to overcome this issue is
to try to solve a more general equation of the form
\[ F(u, \xi) = 0, \]
where $\omega_\epsilon + \ddb u$ is a K\"ahler metric and
$\xi\in\mathfrak{g}$. The operator $F$ is constructed so that if
$F(u,\xi)=0$ and $\xi\in\mathfrak{g}_p$ is in the stabilizer of $p$, 
then $\omega_\epsilon + \ddb u$ is an extremal metric, which has
constant scalar curvature if $\xi=0$. At the
same time the linearization of $F$ is surjective. One can then show
that for sufficiently small $\epsilon$, for every point $p\in M$ we
can find a solution $(u_{\epsilon, p}, \xi_{\epsilon, p})$ of the
corresponding equation. The search for cscK metrics is then reduced to
the finite dimensional problem of finding zeros of the map
\[ \begin{aligned}
  \mu_\epsilon : M &\to \mathfrak{g} \\
                         p &\mapsto \xi_{\epsilon, p}, 
\end{aligned}\]
since if $\mu_\epsilon(p)=0$, then we have found a cscK metric on
$\Bl_p M$ in the class $\Omega_\epsilon$. More generally to find
extremal metrics we need to find $p$ such that
$\mu_\epsilon(p)\in\mathfrak{g}_p$. 

At this point it becomes important to understand better what the map
$\mu_\epsilon$ is, and for this one needs to construct better
approximate solutions than our crude attempt $\omega_\epsilon$
above. In turn this requires more precise expansions of the metrics
$\epsilon^2\eta$ and $\omega$ than what we had in Equations
\eqref{eq:e2eta} and \eqref{eq:omegaexpand}. For the Burns-Simanca
metric, according to Gauduchon~\cite{Gau12} we have
\[
\epsilon^2\eta = \ddb\Big[ |z|^2 - d_0 \epsilon^{2m-2}|z|^{4-2m} + d_1
\epsilon^{2m} |z|^{2-2m} + O(\epsilon^{4m-4}|z|^{6-4m})\Big], 
\]
where $d_0, d_1 > 0$, while for the metric $\omega$ we have
\[ 
\omega = \ddb\Big[ |z|^2 + A_4(z) + A_5(z) + O(|z|^6)\Big], \]
where $A_4$ and $A_5$ are quartic and quintic expressions. Essentially
$A_4$ is the curvature of $\omega$ at $p$, while $A_5$ is its
covariant derivative at $p$. The way to obtain better approximate
solutions than $\omega_\epsilon$ is to preserve more terms in these
expansions rather than multiplying them all with cutoff functions. In
practice this involves modifying the metric $\omega$ on the punctured
manifold $M\setminus\{p\}$ and $\epsilon^2\eta$ on $\Bl_0\mathbf{C}^n$ to
incorporate new terms in their K\"ahler potentials that are asymptotic
to $-d_0\epsilon^{2m-2}|z|^{4-2m} + d_1\epsilon^{2m}|z|^{2-2m}$ and
$A_4(z) + A_5(z)$ respectively. 

The upshot is that we can obtain an expansion for $\mu_\epsilon$ which
is roughly of the form
\[ \label{eq:muepsilon}\mu_\epsilon(p) = \mu(p) + \epsilon^2\Delta\mu(p) +
O(\epsilon^\kappa) \]
for some $\kappa >2$, where $\mu:M\to\mathfrak{g}$ is the moment map
for the action of $G$ on $M$, and $\Delta\mu$ is its
Laplacian. At this point one can exploit the special structure of
moment maps to show that if $\mu(p) + \epsilon^2\Delta\mu(p)$ is in
the stabilizer $\mathfrak{g}_p$, and $\epsilon$ is sufficiently
small, then there is a point $q\in G^c\cdot p$ in the orbit of $p$
under the complexified group such that $\mu_\epsilon(q)\in
\mathfrak{g}_q$. Since $\Bl_p M$ is biholomorphic to $\Bl_q M$ in this
case, we end up with an extremal metric on $\Bl_pM$. Under the
K-polystability assumption this extremal metric is easily seen to have
constant scalar curvature. 

Finally, if $\mu(q) + \epsilon^2\Delta\mu(q) \not\in
\mathfrak{g}_q$ for any $q\in G^c\cdot p$ and sufficiently small
$\epsilon$, then the Kempf-Ness principle~\cite{MFK94}
 relating moment maps to GIT
stability can be exploited to find a $\mathbf{C}^*$-action on $M$
which induces a destabilizing test-configuration for $\Bl_p M$. 

There are several interesting problems which we hope to address
in future work.
\begin{enumerate}
  \item 
 One should extend
Theorem~\ref{thm:stabblowup} to the case when $\dim M=2$ and to
general extremal metric. In principle both of these extensions should
follow from a more refined expansion of the function
$\mu_\epsilon$ than what we have in Equation~\eqref{eq:muepsilon}, 
but it may be more practical to find a different, more
direct approach.
\item Can one obtain
similar existence results for blow-ups along higher dimensional
submanifolds?
\item If $M$ is an arbitrary compact K\"ahler manifold, is it possible
  to construct a cscK metric on the blowup of $M$ in a sufficiently
  large number of points? This would be analogous to Taubes's
  result~\cite{Taubes92} on the existence of anti-self-dual metrics on
  the blowup of a 4-manifold in sufficiently many points. See
  Tipler~\cite{Tip14} for a related result for toric surfaces, where
  iterated blowups are also allowed. 
\end{enumerate}

\section{What if no extremal metric exists?}\label{sec:noextremal}

Even if $M$ does not admit an extremal metric in a class $c_1(L)$, it is natural to try
minimizing the Calabi functional. That this is closely related to the
algebraic geometry of $(M,L)$ is suggested by the following result, 
analogous to a theorem due to Atiyah-Bott~\cite{AB83} in the case of
vector bundles.

\begin{theorem}[Donaldson~\cite{Don05}]\label{thm:Callower}
Given a polarized manifold $(M,L)$, we have
  \[\label{eq:Callower} \inf_{\omega\in c_1(L)} \Vert S(\omega) -
  \underline{S}\Vert_{L^2} \geqslant \sup_{\chi}
  -c_n \frac{\mathrm{Fut}(\chi)}{ \Vert\chi\Vert}, \]
  where the supremum runs over all test-configurations for $(M,L)$ with
  $\Vert \chi\Vert > 0$, and $c_n$ is an explicit dimensional constant.  
\end{theorem}

Donaldson also conjectured that in fact equality holds in
\eqref{eq:Callower}. When $M$ admits an extremal metric
$\omega_{e}\in c_1(L)$,
then it is easy to check that
\[ \Vert S(\omega_{e}) -\underline{S}\Vert_{L^2} = -c_n\frac{
  \mathrm{Fut}(\chi_{e})}{\Vert\chi_{e}\Vert}, \]
where $\chi_{e}$ is the product configuration built from the
$\mathbf{C}^*$-action induced by $\nabla S(\omega_{e})$. In other
words equality holds in \eqref{eq:Callower} in this case. When $(M,L)$
admits no extremal metric, there is little known, except for the case
of a ruled surface~\cite{GSz07_1} where we were able to perform
explicit constructions of metrics and test-configurations to realize
equality in \eqref{eq:Callower}. Note that in the case of vector
bundles the analogous conjecture is known to hold (i.e. equality in
\eqref{eq:Callower}) by Atiyah-Bott~\cite{AB83} over Riemann surfaces,
and Jacob~\cite{J11} in higher dimensions. 

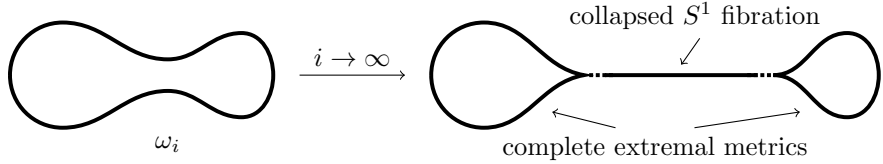
\begin{figure}[h]
\centering
 \begin{tikzpicture}[scale=0.7]
      \draw[line width=1.5pt] (0,0) to [out=90, in=180] (1,1)
      to [out=0, in=180] (3,0.3)
      to [out=0, in=180] (4.4, 0.8)
      to [out=0, in=90] (5,0)
      to [out=270, in=0] (4.4, -0.8)
      to [out=180, in=0] (3,-0.3)
      to [out=180, in=0] (1,-1)
      to [out=180, in=270] (0,0);
      \draw[->] (5.5,0) -- (7.5,0);
      \draw (3, -1) node[below] {$\omega_i$};
      \draw (6.5,0) node[above] {$i\to \infty$};
      \draw[line width=1.5pt] (8,0) to [out=90, in=180] (9,1)
      to [out=0, in=170] (10.9,0)
      to [out=190, in=0] (9,-1)
      to [out=180, in=270] (8,0);
      \draw[line width=1.5pt, style=densely dotted] (11,0) -- (11.5, 0);
      \draw[line width=1.5pt] (11.3,0) -- (14.2,0);
      \draw[line width=1.5pt, style=densely dotted] (14,0) -- (14.5, 0);
      \draw[line width=1.5pt] (14.5,0) to [out=0, in=180] (15.9, 0.8)
      to [out=0, in=90] (16.5,0)
      to [out=270, in=0] (15.9, -0.8)
      to [out=180, in=0] (14.5,0);

      \draw[<-] (12.75, 0.2) to (13.1, 0.7);
      \draw (13,0.7) node[above] {collapsed $S^1$ fibration};
      \draw[<-] (10.3,-0.7) to (11.5, -1);
      \draw[<-] (15, -0.7) to (13, -1);
      \draw (12.25, -1) node[below] {complete extremal metrics};
    \end{tikzpicture}
 \caption{The fiber metrics of a minimizing sequence when $m > k_2$.}
\label{fig:case3metrics}
\end{figure}

To describe our result, let $\Sigma$ be a genus 2 curve, and $M =
\mathbf{P}(\mathcal{O} \oplus \mathcal{O}(1))$, where $\mathcal{O}(1)$
denotes any degree one line bundle over $\Sigma$. For any real number
$m > 0$, we have a K\"ahler class $\Omega_m$ on $M$, defined by 
\[ \Omega_m = [F] + m[S_0], \]
where $[F], [S_0]$ denote Poincar\'e duals to the homology classes of a
fiber $F$ and the zero section $S_0$. Up to scaling we obtain all
K\"ahler classes on $M$ in this way. Depending on the value of $m$, in
\cite{GSz07_1}
we observed 3 qualitatively different behaviors of a minimizing
sequence for the Calabi functional in $\Omega_m$. There are explicitly
computable numbers $0 < k_1 < k_2$ and minimizing sequences $\omega_i$
for the Calabi functional in $\Omega_m$ with the following properties:

\begin{enumerate}
\item When $m < k_1$, then the $\omega_i$ converge to the extremal
  metric in $\Omega_m$ whose existence was shown by
  T\o{}nnesen-Friedman~\cite{TF98}. 
\item When $k_1\leqslant m \leqslant k_2$ then 
  suitable pointed limits of the $\omega_i$
  are complete extremal metrics on $M\setminus S_0$ and
  $M\setminus S_\infty$. Here $S_\infty$ is the infinity section, and
  the volumes of the two complete extremal metrics add up to the
  volumes of the $\omega_i$.
\item When $m > k_2$, then suitable pointed limits of the $\omega_i$
  are either $\Sigma\times \mathbf{R}$, or complete extremal metrics
  on $M\setminus S_0$ or $M\setminus 
  S_\infty$. In the first case a circle
  fibration collapses, and the sum of the volumes of the two complete extremal
  metrics is strictly less than the volume of the
  $\omega_i$. Figure~\ref{fig:case3metrics} illustrates the behavior
  of the metrics $\omega_i$ when restricted to a 
  $\mathbf{P}^1$ fiber. 
\end{enumerate}

We interpret cases 2 and 3 as saying that a minimizing sequence breaks
the manifold into several pieces. Some of the pieces admit complete
extremal metrics, but others display more complicated collapsing
behavior. Having such infinite diameter limits, and possible
collapsing is in stark contrast with the case of Fano manifolds that
we discussed in Section~\ref{sec:existence}. 

The sequences of metrics $\omega_i$ above can be written down
explicitly using the momentum construction developed in detail by
Hwang-Singer~\cite{HS02}. To show that these sequences actually
minimize the Calabi energy, one needs to consider the right hand side
of \eqref{eq:Callower}, and construct corresponding sequences of
test-configurations $\chi_i$ such that
\[\label{eq:limeq}
 \lim_{i\to\infty} \Vert S(\omega_i) - \underline{S}\Vert_{L^2} =
\lim_{i\to\infty}
-c_n\frac{\mathrm{Fut}(\chi_i)}{\Vert\chi_i\Vert}. \]

\begin{figure}[h]
\centering
    \begin{tikzpicture}
      \draw[line width=1.5pt] (0,2) to [out=290, in=80] (-0.2,-2);
      \draw[style=dashed, ->]  (1,0) -- (3,0);
      \draw[line width=1.5pt] (4,1.8) to [out=320, in=110] (5.7,0.3);
      \draw[line width=1.5pt] (5.8,0.6)  to (5.3, 0.1);
      \draw[line width=1.5pt] (5.3,0.3)  to (5.7, -0.1);
      \draw[line width=1.5pt] (5.7,0.1)  to (5.3, -0.3);
      \draw[line width=1.5pt] (5.3,-0.1)  to (5.7, -0.5);
      \draw[line width=1.5pt] (5.7,-0.3)  to (5.3, -0.7);
      \draw[line width=1.5pt] (5.3,-0.5)  to (5.7, -0.9);
      \draw[line width=1.5pt] (5.7,-0.7)  to [out=225, in=70] (4.5,
      -2);
      \draw[<->] (6, 0.6) to (6, -0.8);
      \draw (6.2, 0) node[right, text width=4cm] {more and more
        components as $i\to\infty$};
    \end{tikzpicture}
\caption{The test-configuration $\chi_i$ restricted to a
  $\mathbf{P}^1$ fiber.}
\label{fig:case3degen}
\end{figure}
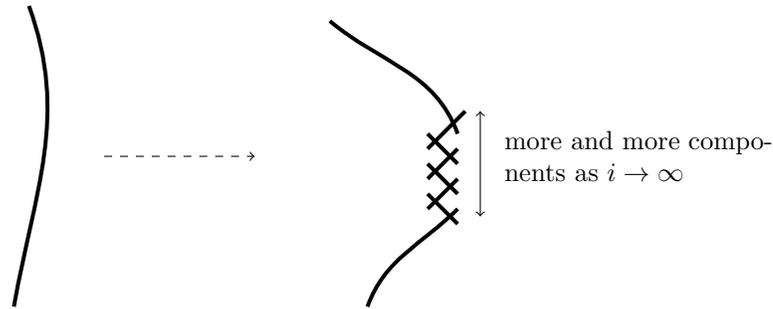

For this to make sense we need to assume that $m$ is rational, so
that a multiple of $\Omega_m$ is an integral class.
Such a sequence $\chi_i$ can be constructed explicitly, and in
the case when $m > k_2$, the
exponents of the test-configurations $\chi_i$ tend to
infinity with $i$. In other words, we need to embed $M$ into larger and larger
projective spaces in order to realize $\chi_i$ as a
degeneration in projective space. The reason is that the central fiber
of $\chi_i$ is a normal crossing divisor consisting of a 
chain of a large number of components isomorphic to
$M$, with the infinity section of each meeting
the zero section of the next one. The number of components goes to
infinity with $i$. Figure~\ref{fig:case3degen} illustrates $\chi_i$
restricted to a $\mathbf{P}^1$-fiber of $M$. 

From Equation~\eqref{eq:limeq} together with
Theorem~\ref{thm:Callower}  we obtain the following.

\begin{theorem}
  For the ruled surface $M$ equality holds in
  Equation~\ref{eq:Callower} for any polarization $L$. 
\end{theorem}

To conclude this section we point out that already in this
example we cannot take a limit of the sequence $\chi_i$ in the space
of test-configurations, because the exponents go to infinity. However
there is a filtration $\chi$, such that $\chi_i$ is the induced
test-configuration of exponent $i$, and in this sense the limit of the
$\chi_i$ exists as a filtration. This filtration achieves the supremum
on the right hand side of \eqref{eq:Callower} and it is natural to ask
whether such a ``maximally destabilizing'' 
filtration always exists. In view of the work of
Bruasse-Teleman~\cite{BT05} this filtration, if it exists, should be
viewed as analogous to the Harder-Narasimhan filtration of an unstable
vector bundle. 
\providecommand{\bysame}{\leavevmode\hbox to3em{\hrulefill}\thinspace}
\providecommand{\MR}{\relax\ifhmode\unskip\space\fi MR }
\providecommand{\MRhref}[2]{%
  \href{http://www.ams.org/mathscinet-getitem?mr=#1}{#2}
}
\providecommand{\href}[2]{#2}

\end{document}